\documentclass{ejm}

\pdfoutput=1
\usepackage{graphicx}
\usepackage{amsmath}
\usepackage{float}
\usepackage{hyperref}
\usepackage{amsfonts}

\title[Deep Learning for Global Coordinate Transformations that Linearize PDEs]{Deep Learning Models for Global Coordinate Transformations that Linearize PDEs}
\author[Gin, Lusch, Brunton, \and Kutz]{Craig Gin$\,^1$, Bethany Lusch$\,^2$, Steven L. Brunton$\,^{1,3}$, \and J. Nathan Kutz$\,^1$}

\affiliation{%
  $^1\,$Department of Applied Mathematics, University of Washington, Seattle, WA, 98195, USA\\
  $^2\,$Argonne Leadership Computing Facility, Argonne National Laboratory, Lemont, IL, USA\\
  $^3\,$Department of Mechanical Engineering, University of Washington, Seattle, WA, 98195, USA}

\begin{document}
\maketitle
\thispagestyle{empty}
\begin{abstract}
We develop a deep autoencoder architecture that can be used to find a coordinate transformation which turns a nonlinear PDE into a linear PDE.  Our architecture is motivated by the linearizing transformations provided by the Cole-Hopf transform for Burgers equation and the inverse scattering transform for completely integrable PDEs.  By leveraging a residual network architecture, a near-identity transformation can be exploited to encode intrinsic coordinates in which the dynamics are linear. The resulting dynamics are given by a Koopman operator matrix $\mathbf{K}$. The decoder allows us to transform back to the original coordinates as well. Multiple time step prediction can be performed by repeated multiplication by the matrix $\mathbf{K}$ in the intrinsic coordinates.  We demonstrate our method on a number of examples, including the heat equation and Burgers equation, as well as the substantially more challenging Kuramoto-Sivashinsky equation, showing that our method provides a robust architecture for discovering interpretable, linearizing transforms for nonlinear PDEs.
\end{abstract}


\begin{subjclass}[2010]
35A22, 35A35, 37M99, 65P99, 68T99
\end{subjclass}

\section{Introduction}

 \begin{figure}[t]
 \centering
\includegraphics{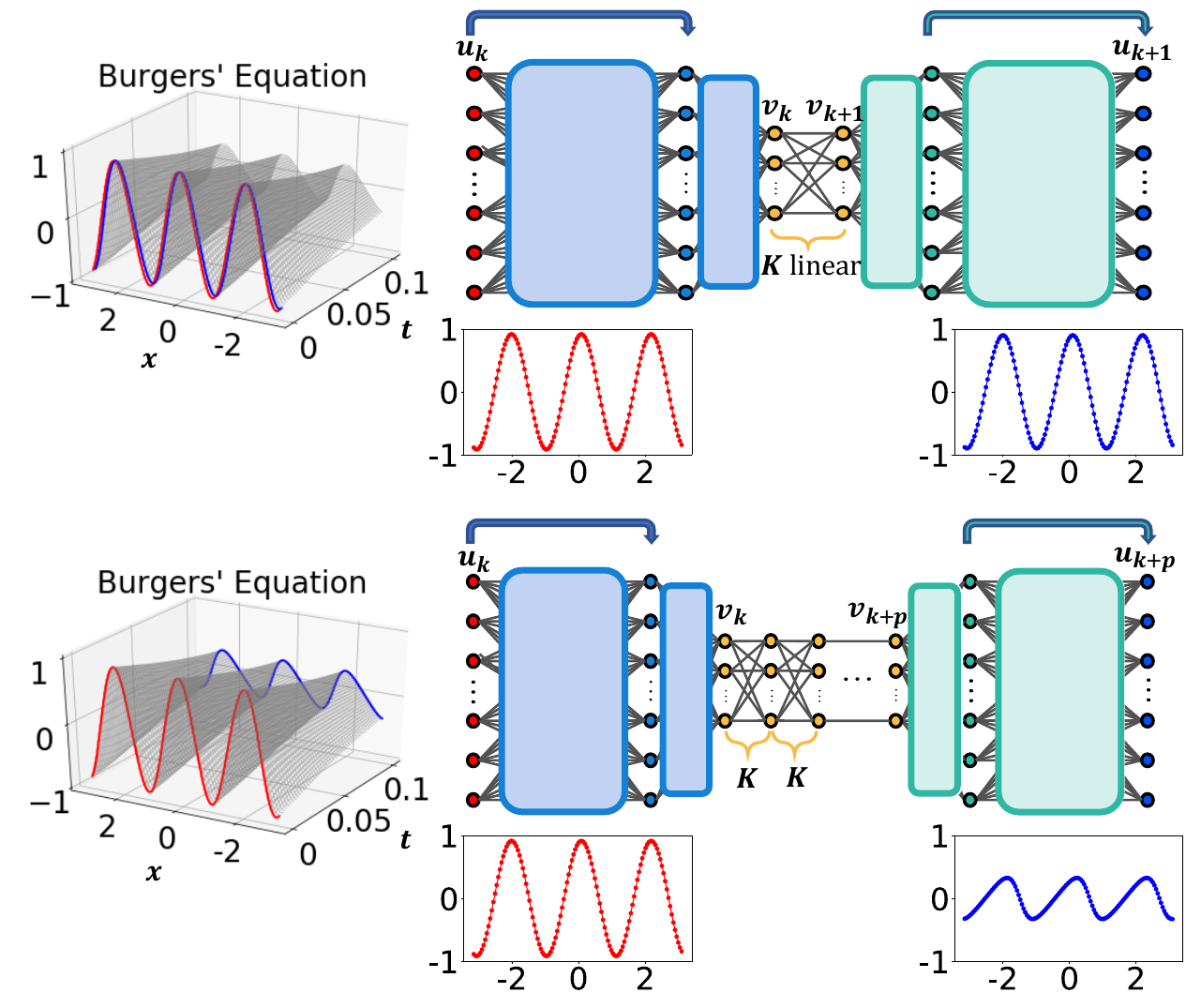}
 \caption{A deep autoencoder is used to find coordinate transformations to linearize PDEs. The encoder finds a set of intrinsic coordinates for which the dynamics are linear. Then the dynamics are given by a matrix $\mathbf{K}$. The decoder transforms back to the original coordinates. Multiple time step prediction can be performed by repeated multiplication by the matrix $\mathbf{K}$ in the intrinsic coordinates.}
\end{figure}
 

Partial differential equations (PDEs) provide a theoretical framework for modeling spatio-temporal systems across the biological, physical and engineering sciences.  Analytic solution techniques are readily available for PDEs that are linear and have constant coefficients~\cite{haberman}.  These PDEs include canonical models such as the heat equation, wave equation and Laplace's equation which are amenable to standard separation of variable techniques and linear superposition.  In contrast, there is no general mathematical architecture for solving nonlinear PDEs as methods like separation of variables fail to hold, thus recourse to computational solutions is necessary.  There are a few, but notable, exceptions:  (i) the Cole-Hopf transformation~\cite{hopf50,cole51} for solving diffusively regularized Burgers equation, and (ii) the {\em Inverse Scattering Transform} (IST)~\cite{ist} for solving a class of completely integrable PDEs such as Korteweg deVries (KdV), nonlinear Schr\"odinger (NLS), Klein-Gordon, etc.  The success of Cole-Hopf and IST is achieved by providing a {\em linearizing} transformation of the governing nonlinear equations.  Thus these methods provide a transformation to a new coordinate system where the dynamics is characterized by a linear PDE, and for which recourse can be made to the well-established methods for solving linear PDEs.  Recent efforts have shown that it is possible to use neural networks to discover advantageous coordinate transformations for dynamical systems~\cite{lusch2018deep,Champion2019,Wehmeyer2017arxiv,Mardt2017arxiv,Takeishi2017nips,Yeung2017arxiv,Otto2017arxiv,Li2017chaos,dsilva2018parsimonious}.  Indeed, such architectures provide a strategy to discover linearizing transformations such as the Cole-Hopf and IST.   In this manuscript, we develop a principled method for using neural networks for discovering coordinate transformations for linearizing nonlinear PDEs, demonstrating that transformations such as Cole-Hopf can be discovered from data alone with proper architecting of the network.

Linearizing transforms fall broadly under the aegis of Koopman operator theory~\cite{Koopman1931pnas} which has a modern interpretation in terms of dynamical systems theory~\cite{Mezic2004,Mezic2005nd,budivsic2012applied,Mezic2013arfm}.   Koopman operators can only be constructed explicitly in limited cases~\cite{brunton_koopman_2016}, however
{\em Dynamic Mode Decomposition} (DMD)~\cite{Schmid2010jfm} provides a numerical algorithm to provide an approximation to the Koopman operator~\cite{Rowley2009jfm}, with many recent extensions improving on the DMD approximation~\cite{Kutz2016book}.   Importantly, many of the advantageous transformations highlighted above attempt to construct Koopman embeddings for the dynamics using neural networks~\cite{lusch2018deep,Wehmeyer2017arxiv,Mardt2017arxiv,Takeishi2017nips,Yeung2017arxiv,Otto2017arxiv,Li2017chaos}.  This is in addition to enriching the observables of DMD~\cite{noe2013variational,nuske2014variational,Williams2015jnls,Williams2014arxivA,klus2017data,kutzPDE,page2018koopman}.
Thus neural networks have emerged as a highly effective mathematical architecture for approximating complex data~\cite{ml3,GoodfellowDL}.  Its universal approximation properties~\cite{Cybenko.1989,Hornik.1990} are ideal for learning the coordinate transformations required for linearizing nonlinear PDEs.  NNs have also been used in this context to discover time-stepping algorithms for complex systems~\cite{rico1995nonlinear,gonzalez1998identification}, which is a slightly different, but related task to what is advocated here.  Critical to the success of a neural network is imposing proper constraints, or regularizers, for the desired task.  For spatio-temporal systems modeled by PDEs, there are a number of constraints that are required for success.  These constraints are largely motivated by domain knowledge, thus producing a physics-informed machine learning architecture for PDEs.   Specifically, we identify four critical components for successfully training a neural network for nonlinear PDEs:  (i) The appropriate neural network architecture for the desired task, (ii) A method for handling the identity transformation as ideas of near-identity transformation are critical, (iii) Domain knowledge constraints for the PDEs, and (iv) Judiciously chosen spatio-temporal data for training the network.

Imposing structure on the NN is the first critical design task.   As shown in recent works, autoenconders (AEs) are a physics informed choice as they have been shown to be able to take input data from the original high-dimensional input space to the new coordinates which are at the intrinsic rank of the underlying dynamics~\cite{lusch2018deep,Champion2019,pan2019}.     AEs allow for the computation of a nonlinear dimensionality reduction in constructing a linear reduced order model.   The AE mapping is significantly improved if the identity map can be easily implemented.  The success of near-identity transformations for transforming dynamical systems motivates the leading order identity mapping~\cite{neu1980method,wiggins2003introduction}.   As will be shown, nonlinear activation functions such as ReLU make it difficult to produce the identity operator.  The success of {\em Deep Residual Networks} (DRN)~\cite{he2016deep} motivate our approach to handling the identity map.  We have found that the AE architecture which leverages the DRN structure provides a physics-informed architecture that leverages the concepts of near-identity transforms along with the low-dimensional intrinsic rank of the physics itself.

Finally the data must be judiciously chosen.  As stated by S. Mallat:  {\em Supervised learning is a high-dimensional interpolation problem}~\cite{mallat2016understanding}.  Thus a wide range of initial conditions and trajectories must be sampled in order to construct a robust and accurate linearizing transform.  The success of training the AE, and indeed of any neural network architecture, hinges on training with a broad enough class of data so that new  spatio-temporal trajectories can leverage the interpolation powers of the AE infrastructure.  We emphasize again that NNs are interpolation engines and fail in extrapolatory scenarios.

We demonstrate our method on a set of prototype PDE models:  the heat equation, Burgers equation and the Kuramoto-Shivasinsky (KS) equation.   For the first two, we have ground truth analytic results to validate the methodology.  For the KS equation, we provide a first mapping of the PDE to a linearized system which can be completely inverted, although a number of Galerkin schemes have been advocated to characterize the KS dynamics~\cite{foias1994some}.  These examples demonstrate the power of the methodology and the required regularizations necessary for training NNs for learning coordinate embeddings.

\section{Numerical Timestepping and Dynamics}

Given the spatio-temporal nature of data generated by PDEs, we develop neural network architectures that exploit various aspect of numerical time-stepping schemes.  This begins with understanding the near-identity transformation generated from numerical discretization schemes for solving PDEs.  

\subsection{Identity Functions and Neural Networks}
\begin{figure}[t]
 \centering
\includegraphics{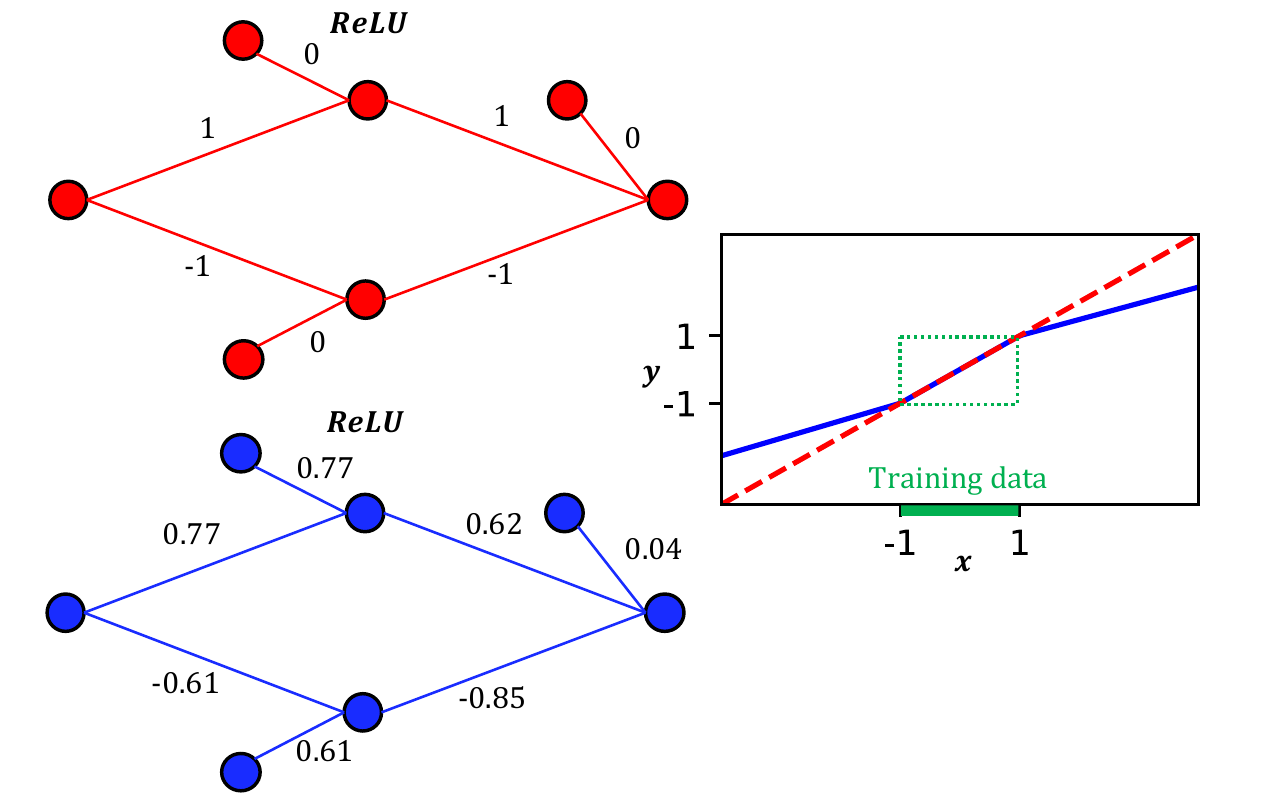}
 \caption{We train a seven-parameter nonlinear neural network to fit $f(x)=x$ with data from $[-1,1]$. Although there is an optimal solution (see the red network) that is valid for all $x \in \mathbb{R}$, the training algorithm chooses parameters such as those labeled on the blue network. On the right, the function learned by the blue network is plotted in blue. This is compared to the identity function, plotted in red. The training domain is annotated in green. We note that the trained network should not be used for extrapolation.} \label{Fig2}
\end{figure}

\begin{figure}[t]
 \centering
\includegraphics{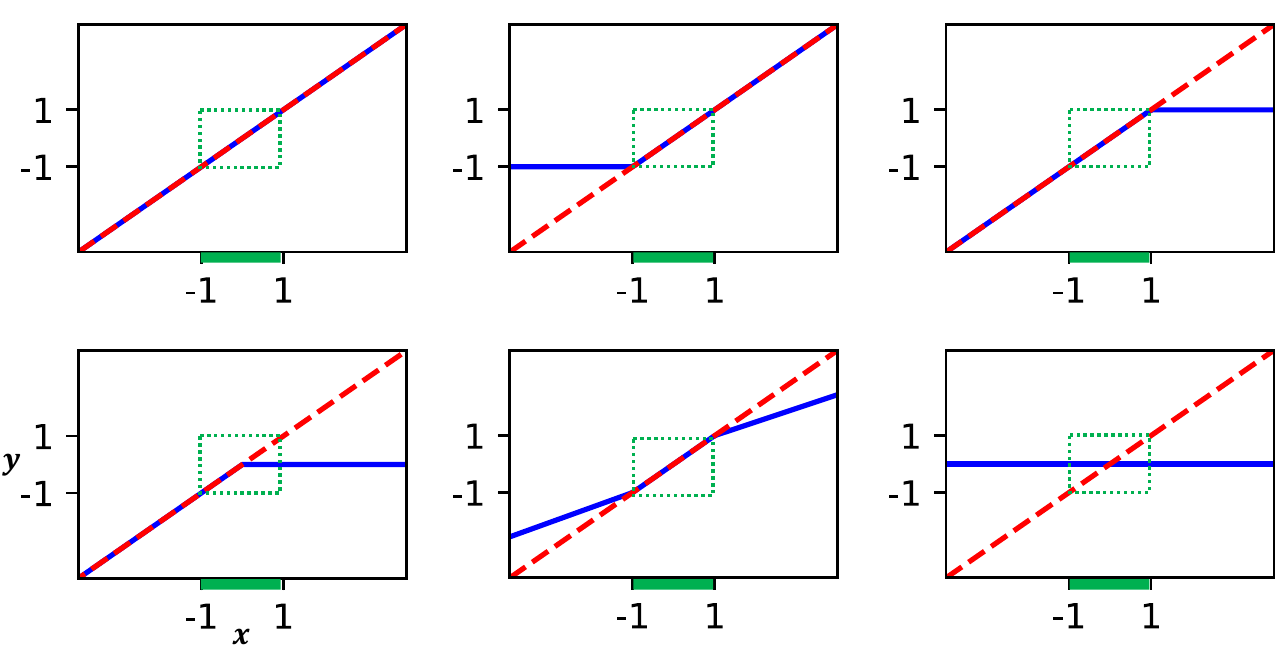}
 \caption{We train a small neural network (shown in Figure \ref{Fig2}) to fit $f(x)=x$ with data from $[-1,1]$. We show results from six training trials with different random initializations. In each of the six plots, the function learned by the network is plotted in blue. This is compared to the identity function, plotted in red. The training domain is annotated in green. We note that the trained network should not be used for extrapolation, and that sometimes the trained network is not accurate even on the training data.} \label{Fig3}
\end{figure}

Although training a neural network to represent $f(x) = x$ may seem simple, we show that nonlinear activation functions such as ReLU make this non-trivial. We consider a one-dimensional dataset $x$ composed of random real numbers from $[-1,1]$. We assign each point a label $y = x$. Then we create a small seven-parameter neural network with a one-node input later, two-node hidden layer with a ReLU activation, and a one-node linear output layer. Each node in the hidden layer and output layer have a bias term. This network is depicted in Figure \ref{Fig2}. The red network in Figure \ref{Fig2} is labeled with parameters that would map input $x$ to output $y=x$ perfectly for any $x \in \mathbb{R}$. However, when we train this network on our training data from $[-1,1]$, different parameters are chosen, such as those shown on the blue network in Figure \ref{Fig2}. This trained network is accurate in the training domain of $[-1,1]$ but diverges from $f(x)=x$ outside of that domain. Further, as demonstrated in Figure \ref{Fig3}, if we train the network repeatedly, with different random initializations of the parameters, we obtain different results. In some cases, the training loss plateaus and the training fails to accurately fit the training data at all. We have two observations from this experiment:
\begin{enumerate}
\item We are reminded that neural networks cannot be relied upon to extrapolate outside of the training domain. 
\item Even though it is possible to represent the identity function with this nonlinear neural network, it is non-trivial for the training algorithm to fit the data.
\end{enumerate}
This latter observation contributes to our motivation to use a residual neural network structure. Although in theory a near-identity function $g(x)$ may be representable by a neural network, the optimization may be more successful at learning the residual $h(x)$ where $g(x) = x + h(x)$.

\pagebreak
\subsection{Building Networks for Time-stepping Dynamics}
In what follows, we consider several PDEs, typically nonlinear, of the form
\begin{equation}\label{continuousPDE}
u_t = F(u,u_x,u_{xx}, \cdots, x,t).
\end{equation}
where $F(\cdot)$ characterizes the governing equations.
When discretized in space and time, the PDE becomes a finite-dimensional, discrete-time dynamical system of the form~\cite{Kutz:2013}
\begin{equation}\label{discreteDynSys}
 \mathbf{u_{k+1}} = \mathbf{F}(\mathbf{u_k}),
\end{equation}
where $\mathbf{u_k}$ is a spatial discretization of the function $u(x,t)$ at time $t_k$ where $t_k = k \Delta t$ for some time-step $\Delta t$. In the case that the function $\mathbf{F}$ is linear, the future values of the state $\mathbf{u}$ can be solved for exactly using a spectral expansion. However, $\mathbf{F}$ is typically nonlinear, and there is no general framework for solving nonlinear systems. 

Koopman operator theory provides a means to linearize a nonlinear dynamical system. The Koopman operator, $\mathcal{K}$, is a linear operator that acts on the infinite-dimensional space of observables $\mathbf{v} = g(\mathbf{u})$ of the state. The Koopman operator is defined by
\begin{equation}
 \mathcal{K}g(\mathbf{u_k}) = g(\mathbf{u_{k+1}}).
\end{equation}
Therefore, the Koopman operator advances the observables in time. Because of its linearity, the behavior of the Koopman operator is completely determined by its eigenvalues and eigenfunctions. We use deep learning in order to approximate the Koopman eigenfunctions, which satisfy
\begin{equation}
 \varphi(\mathbf{u_{k+1}}) = \mathcal{K}\varphi(\mathbf{u_k}) = \lambda \varphi (\mathbf{u_{k}}).
\end{equation}
Although the Koopman operator acts on an infinite-dimensional space, we can obtain a finite-dimensional approximation by considering the space spanned by finitely many Koopman eigenfunctions. Acting on this space, the Koopman operator is just a matrix. Therefore, Koopman operator theory provides an approach to find an intrinsic coordinate system in which the dynamical system has linear dynamics.

\begin{figure}[t]
\centering
\includegraphics{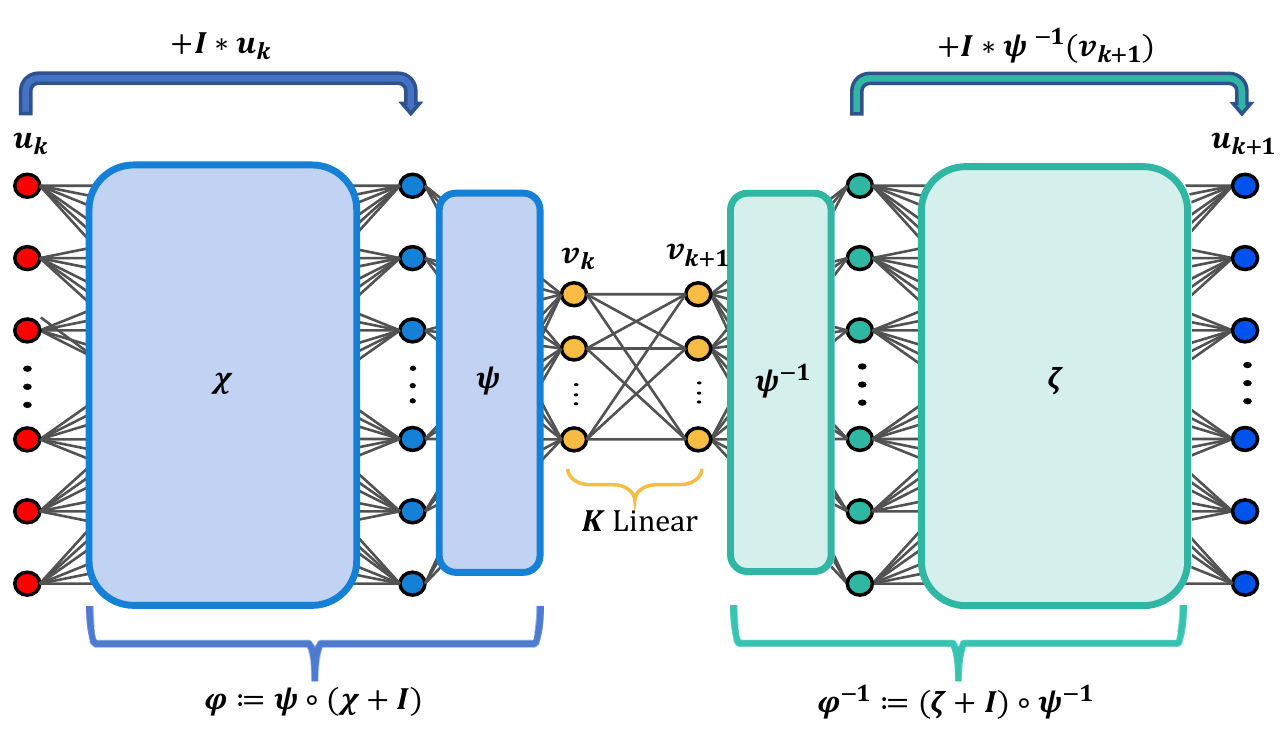}
 \caption{The network architecture used to find the Koopman eigenfunctions. The network consists of an outer encoder, inner encoder, dynamics matrix, inner decoder and outer decoder. The outer encoder and decoder use a residual neural network architecture.}
 \label{fig:NetworkArch}
\end{figure}

In this work, we use the universal approximation properties of neural networks to find such linearizing coordinate transformations. The network architecture that we use is shown in Figure \ref{fig:NetworkArch}.
The input of the the network $\mathbf{u_k}$ is the state vector at time $t_k$ and the output is the state vector at time $t_{k+1}$. The network consists of three parts: (i) the encoder $\mathbf{\varphi}$, (ii) the linear dynamics $\mathbf{K}$, and (iii) the decoder $\mathbf{\varphi}^{-1}$. Both the encoder and decoder are split into two parts. The encoder consists of the outer encoder $\mathbf{\chi} + \mathbf{I}$ and the inner encoder $\mathbf{\psi}$. The outer encoder performs a coordinate transformation into a space in which the dynamics are linear. The inner encoder diagonalizes the system and/or reduces the dimensionality. The inner decoder $\mathbf{\psi}^{-1}$ and the outer decoder $\mathbf{\zeta} + \mathbf{I}$ are the inverses of the inner and outer encoder, respectively.

For a sufficiently small time step, the function $\mathbf{F}(\mathbf{u})$ in equation \eqref{discreteDynSys} is a near-identity function, which is a simple representation of an integration scheme~\cite{rico1995nonlinear,gonzalez1998identification,Kutz:2013,pan2019}. As demonstrated in the previous section, the identity transformation can be difficult to learn using nonlinear neural networks. Therefore, special care must be given to account for the fact that the network is approximating a near-identity function~\cite{neu1980method,wiggins2003introduction}. This is handled in two ways. First, a residual neural network architecture is used for both the outer encoder and outer decoder~\cite{he2016deep}. The coordinate transformation given by the outer encoder is represented as the identity transformation plus a residual which is given by the neural network $\chi$. Therefore, the outer decoder is $\mathbf{\chi} + \mathbf{I}$. Similarly the outer decoder is $\mathbf{\zeta} + \mathbf{I}$ where $\mathbf{\zeta}$ is a neural network. Second, the interior layers of the network, $\mathbf{\psi}$, $\mathbf{K}$, and $\mathbf{\psi}^{-1}$, are all initialized as identity matrices. 

\begin{figure}[t]
 \centering
\includegraphics{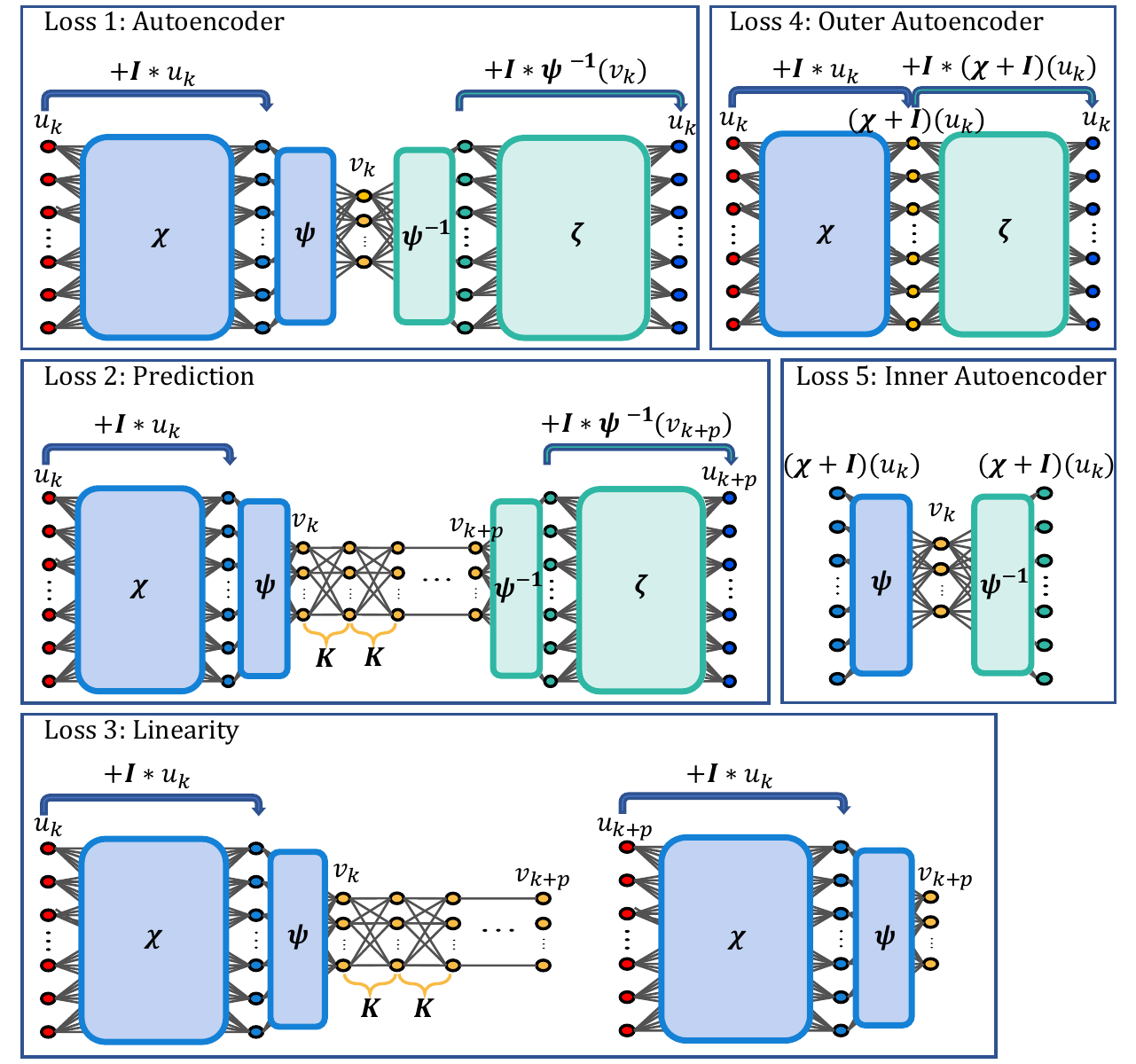}
 \caption{A depiction of the loss functions used for training the neural network.}
 \label{fig:Lossfns}
\end{figure}

The loss function used to train the network is the sum of five different losses. They are depicted in Figure \ref{fig:Lossfns}. Each loss enforces a desired condition:
\begin{enumerate}
 \item {\bf Loss 1: autoencoder loss}. We want an invertible transformation between the state space and intrinsic coordinates for which the dynamics are linear. The transformation into the intrinsic coordinates is given by the encoder $\varphi$ and the transformation back into the state space is given by the decoder $\varphi^{-1}$. Therefore, we wish for the autoencoder $\varphi^{-1} \circ \varphi$ to reconstruct the inputs of the network as closely as possible. This loss is given by $ \left\lVert \mathbf{u_k} - \varphi^{-1}(\varphi(\mathbf{u_k}))\right\rVert$ where $\left\lVert \cdot \right\rVert$ is the mean-squared error averaged over all trajectories in the data. We also add $\ell^2$ regularization.
 \item {\bf Loss 2: prediction loss}. The output of the network should accurately predict the state $\mathbf{u_{k+1}}$ when given the state at the previous time $\mathbf{u_k}$. The loss is given by $\left\lVert \mathbf{u_{k+1}} - \varphi^{-1}(\mathbf{K}\varphi(\mathbf{u_k}))\right\rVert$. Furthermore, we would like to be able to predict multiple time steps into the future by iteratively multiplying by the matrix $\mathbf{K}$. Therefore, in general we have $\left\lVert \mathbf{u_{k+p}} - \varphi^{-1}(\mathbf{K^p}\varphi(\mathbf{u_k}))\right\rVert$. Note that multi-step prediction is done by evolving multiple steps in the intrinsic coordinates and not by passing in and out of the intrinsic coordinates at each time step.
 \item {\bf Loss 3: linearity loss}. The dynamics on the intrinsic coordinates should be linear. Therefore, we enforce a prediction loss within these coordinates: \\$\left\lVert \varphi(\mathbf{u_{k+p}}) - \mathbf{K^p}\varphi(\mathbf{u_k})\right\rVert$.
 \item {\bf Loss 4: outer autoencoder loss}. We wish for the outer encoder and the outer decoder to form an autoencoder. This separates the transformation into the intrinsic coordinates from dimensionality reduction and/or diagonalization which are performed by the inner encoder/decoder. Disentangling these two processes allows for better interpretability in the case of Burgers' equation in Section \ref{sec:Burgers}. The outer autoencoder loss is given by $ \left\lVert \mathbf{u_k} - (\zeta+\mathbf{I})((\chi+\mathbf{I})(\mathbf{u_k}))\right\rVert$.
 \item {\bf Loss 5: outer autoencoder loss}. It is given by \\$\left\lVert (\chi+\mathbf{I})\mathbf{u_k} - \psi^{-1}(\psi((\chi+\mathbf{I})(\mathbf{u_k})))\right\rVert$.
\end{enumerate}

The data for training the neural networks is created by performing numerical simulations of the given PDE. The initial conditions used and discretization details are described for each example below. In each case, we choose the data to be sufficiently diverse as to learn a global coordinate transformation that holds for a wide variety of initial conditions.
Note that our approach is completely data driven - no knowledge of the underlying equations is needed. Therefore, it can be used for experimental data for which the governing equations are unknown.

\section{Heat Equation}

The first PDE that we consider is the one-dimensional heat equation:
\begin{equation}
 u_{t} = u_{xx}, \qquad x \in (-\pi,\pi),
\end{equation}
with periodic boundary conditions. We discretize the spatial domain using $n=128$ equally spaced points.

\begin{figure}[t]
\centering
\includegraphics[width=0.8\textwidth]{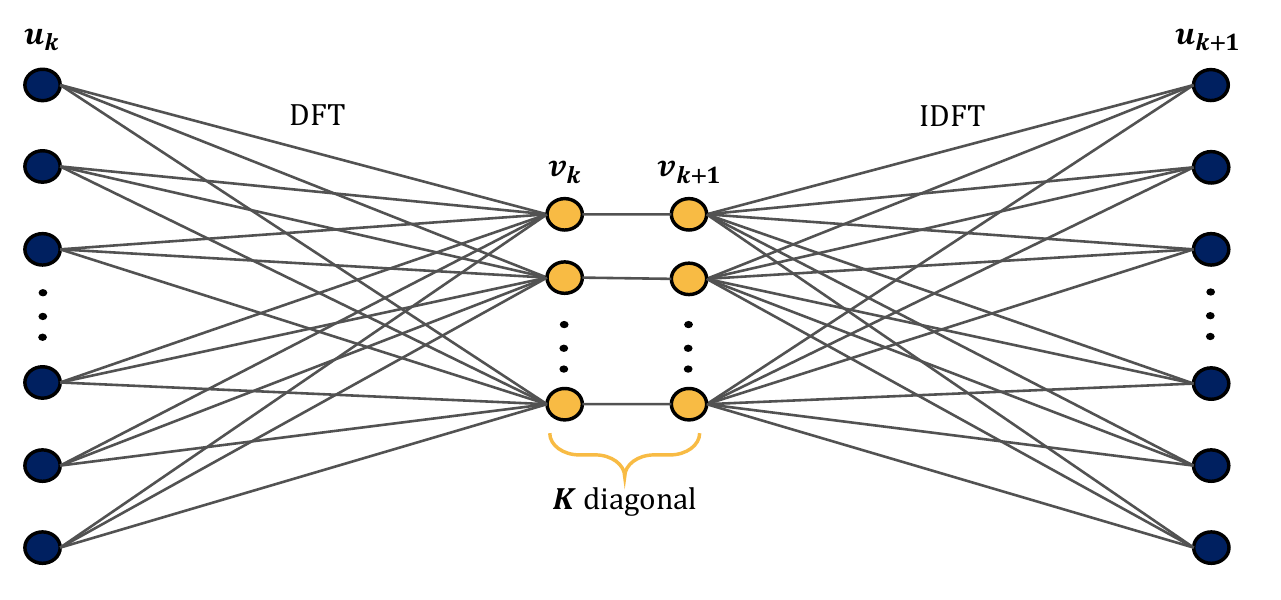}
 \caption{The network architecture for the heat equation. The input and output layers have $n=128$ neurons and the two hidden layers have 21 neurons. The network has no activation functions.}
 \label{fig:HeatArch}
\end{figure}

Because the heat equation is linear, finding the Koopman eigenfunctions amounts to diagonalizing the heat equation. The network architecture that we use for the heat equation is shown in Figure \ref{fig:HeatArch}. There is no outer encoder or decoder because there is no need for a linearizing transformation. Therefore, the entire network is linear and consists of two hidden layers. The widths of the input and output layers match the spatial discretization: $n=128$ neurons in each. The width of the hidden layers is a parameter $r$ that can be adjusted in order to obtain a reduced order model of rank $r$ for the heat equation. We chose $r = 21$. The heat equation is linearized by the Fourier transform so the encoder should mimic a discrete Fourier transform (DFT) that is truncated to include the $r$ most dominant modes. The decoder plays the role of the inverse discrete Fourier transform. The matrix $\mathbf{K}$ that represents the dynamics is forced to be diagonal.

The training data consists of 8000 trajectories from the heat equation. The initial condition in all cases is randomly created white noise, and therefore has no particular structure. The trajectories consist of 50 equally spaced time steps with $\Delta t = 0.0025$. The validation data has the same structure but contains 2000 trajectories.

\begin{figure}[t]
\centering
\includegraphics[width=0.8\textwidth]{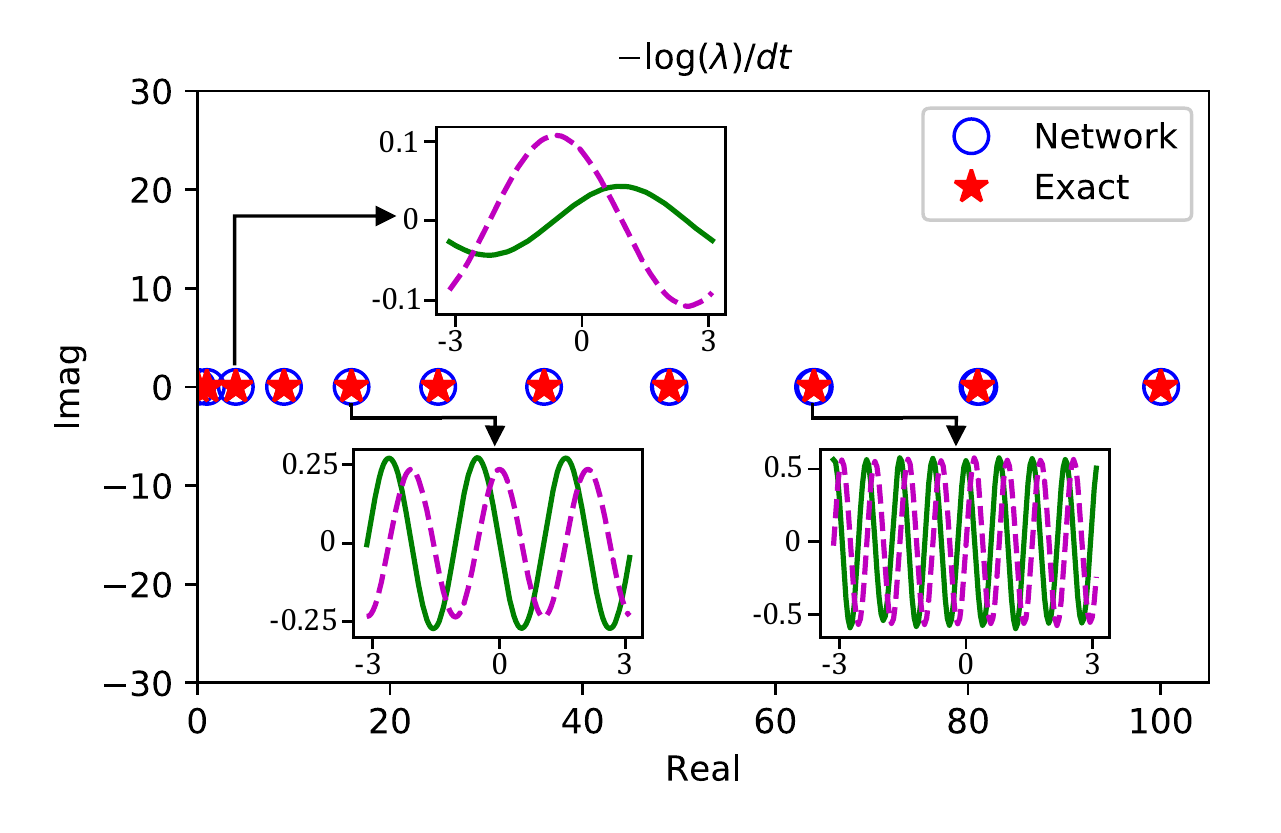}
 \caption{The eigenvalues of the matrix $\mathbf{K}$ from the neural network are plotted along with the exact, discrete-time eigenvalues of the heat equation. The eigenfunctions corresponding to several eigenvalues are shown.}
 \label{fig:HeatEig}
\end{figure}

\begin{figure}[t]
\centering
\includegraphics[width=0.6\textwidth]{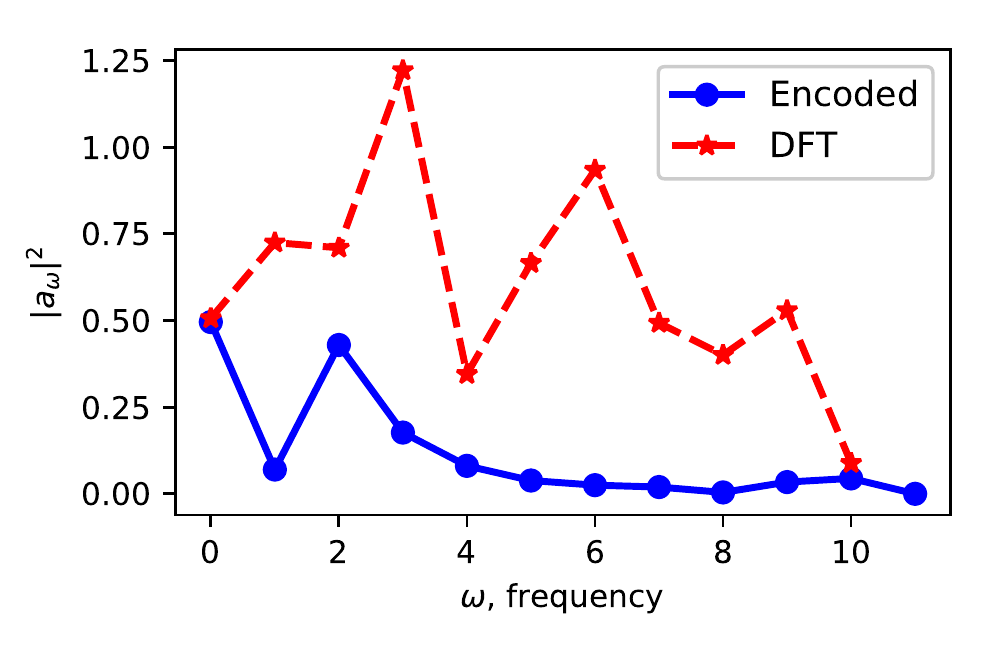}
 \caption{For a particular function, the DFT and the encoded function are shown. The encoder is an alternate transformation to diagonalize the heat equation.}
  \label{fig:HeatDFT}
\end{figure}

For the heat equation, the discrete-time eigenvalues are
\begin{equation}
 \lambda = e^{-\omega^2 \Delta t}, \qquad \omega = 0, \pm 1, \pm 2, ...
\end{equation}
Because the spatial discretization is $n=128$ points, $\omega$ ranges from $-64$ to $63$. The eigenfunctions are $\sin(\omega x)$ and $\cos(\omega x)$ for each positive value of $\omega$. Because the high-frequency waves decay faster than the low-frequency waves, we expect the $21 \times 21$ matrix $K$ to have the eigenvalues 
\begin{equation}
 \lambda = e^{-\omega^2 \Delta t}, \qquad \omega = 0, \pm 1, \pm 2, ..., \pm 10,
\end{equation}
and therefore the eigenvalues satisfy
\begin{equation}\label{logEig}
 -\log(\lambda)/\Delta t = \omega^2, \qquad \omega = 0, \pm 1, \pm 2, ..., \pm 10.
\end{equation}
The eigenvalues of the network are shown in Figure \ref{fig:HeatEig} along with the exact analytical values given by equation \eqref{logEig}. Note that for $\omega \neq 0$, the eigenvalues that are plotted occur in pairs. There is very good agreement between the eigenvalues from the neural network and the exact values. Also in Figure \ref{fig:HeatEig} are the eigenfunctions corresponding to $\omega = \pm 1, \pm 3, \pm 8$. These were plotted by taking the eigenvectors of $\mathbf{K}$ and feeding them through the decoder. As expected, the eigenfunctions are oscillatory and have a frequency that corresponds to the eigenvalues they are associated with. However, note that the eigenfunctions do not all have the same amplitude. This is because the DFT is not a unique transformation to diagonalize the heat equation.  In addition to performing a DFT, the encoder of the network can also scale each of the Fourier modes by any constant. The decoder must then invert the scalings performed by the encoder. The fact that the encoder differs from the DFT is shown in Figure \ref{fig:HeatDFT} which compares the result of feeding a particular function through the encoder versus the DFT of that function. 

\begin{figure}[t]
\centering
\includegraphics{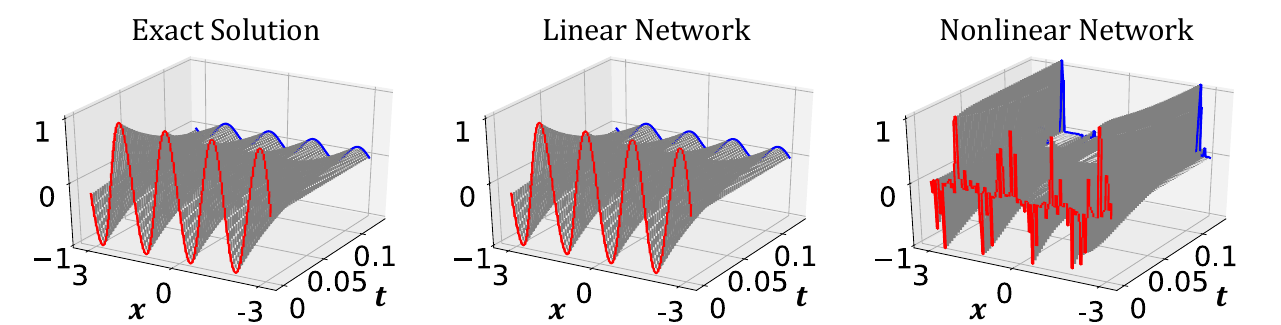}
 \caption{Plots of (a) an exact solution to the heat equation, (b) the prediction given by the linear neural network, and (c) the prediction given by a nonlinear neural network.}
 \label{fig:HeatPred}
\end{figure}

Because the entire network is linear, having the correct eigenvalues and eigenvectors is enough to have a global representation of the heat equation. Therefore, it can be used for prediction for initial conditions that are not represented in the training data. As an example, the network was used to do prediction for the initial condition $u_0(x) = \sin(4x)$. The exact solution and the output of the network are shown in Figure \ref{fig:HeatPred} as the two leftmost plots. The two are indistinguishable. Note that for all network predictions that follow, only the initial condition is needed. To get a prediction at time $t_p = p \Delta t$, the initial condition is encoded, multiplied by the matrix $\mathbf{K}^p$, and then decoded.

It is very important that we leveraged the fact that the heat equation is linear and used a linear neural network. This allowed for the ability to generalize outside of the training data. For contrast, we also trained a nonlinear network by adding two fully-connected layers, the first of which has a ReLU activation function, to both the encoder and decoder. The total validation loss for the nonlinear network $0.0449$ which is less than the $0.0711$ validation loss for the linear network. However, prediction using the nonlinear network for data with a different structure than the training data is very poor, as evidenced by the rightmost plot in Figure \ref{fig:HeatPred}.

Although the heat equation is already linear, it serves as a useful first example because the network was able to find the eigenvalues and eigenfunctions of the heat equation just from the data. In addition, it was able to identify the optimal reduced-order model of rank $r = 21$.

\section{Burgers' Equation}\label{sec:Burgers}

For our second example, we consider the nonlinear PDE known as the diffusively regularized Burgers' equation which can be written as 
\begin{equation}
 u_t + \epsilon u u_x = \mu u_{xx}, \qquad x \in (-\pi,\pi),
\end{equation}
where the parameter $\epsilon$ is the strength of advection and the parameter $\mu$ is the strength of diffusion.  For our results, we use $\epsilon = 10$ and $\mu = 1$. We again use periodic boundary conditions and discretize the spatial domain with $n=128$ equally spaced points. Burgers' equation can be linearized through the Cole-Hopf transformation.  Namely, let $u(x,t)$ be a solution to Burgers' equation and define the function
\begin{equation}\label{CHencode}
 v(x,t) = exp\left[-\frac{\epsilon}{2\mu} \int_0^x u(s,t) ds\right].
\end{equation}
If $u(x,t)$ satisfies Burgers' equation then $v(x,t)$ solves the heat equation. The function $u(x,t)$ can be recovered from $v(x,t)$ by inverting the transformation:
\begin{equation}\label{CHdecode}
 u = -2 \frac{\mu}{\epsilon} \frac{v_x}{v}.
\end{equation}

The neural network architecture used for Burgers' equation is shown in Figure \ref{fig:NetworkArch}. The interpretation of the network in the context of the Cole-Hopf transformation is as follows. The input to the network $\mathbf{u_k}$ is the discretized solution to Burgers' equation at time $t_k$. The outer encoder $\mathbf{\chi} + \mathbf{I}$ performs a linearizing transformation analogous to equation \eqref{CHencode}. The inner encoder $\mathbf{\psi}$ diagonalizes the system and potentially reduces the dimensionality of the system. Since the heat equation governs the dynamics of the intrinsic coordinates, $\mathbf{\psi}$ plays the role of the DFT just like the encoder in the previous section. The matrix $\mathbf{K}$ moves the solution forward in time by one time step. The inner decoder attempts to invert the inner encoder. Finally, the outer decoder $\mathbf{\zeta} + \mathbf{I}$ performs a transformation back to the original coordinates, analogous to equation \eqref{CHdecode}. The output $\mathbf{u_{k+1}}$ is the solution to Burgers' equation at time $t_{k+1} = t_k + \Delta t$.

\begin{figure}[t]
\centering
\includegraphics{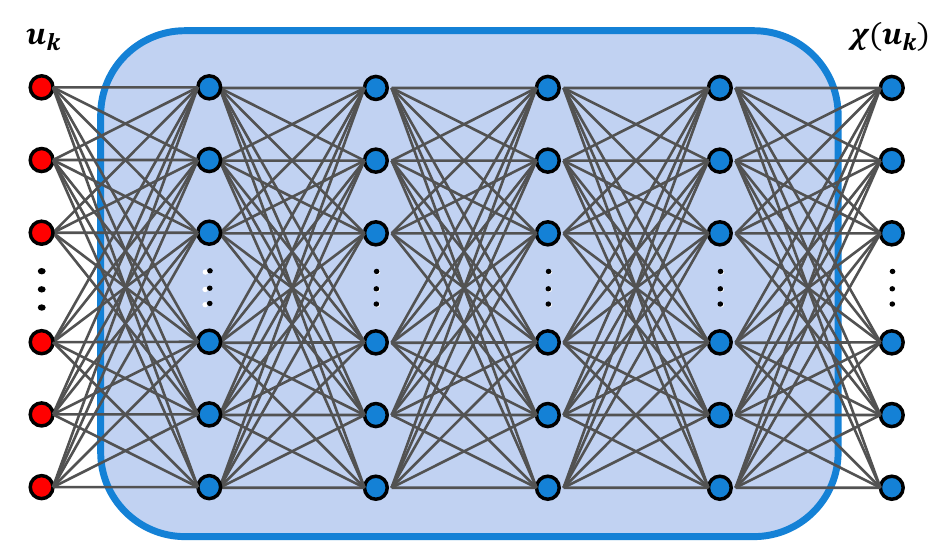}
 \caption{The network architecture for the outer encoder and decoder for Burgers' equation. All layers have 128 neurons. ReLU activation functions are used on all but the last layer.}
 \label{fig:BurgersArch}
\end{figure}

For the inner encoder $\psi$ and inner decoder $\psi^{-1}$, we use a single fully-connected linear layer. For some results presented below, the layer $\mathbf{v_k}$ has the same width as the input and output layers and hence there is no dimensionality reduction. For other results, the dimensionality is reduced by the inner encoder to some lower rank $r$.  
Figure \ref{fig:BurgersArch} shows the network architecture used for $\chi$ in the outer encoder and $\zeta$ in the outer decoder. These parts of the network consist of six fully-connected layers that each contain $n=128$ neurons. All but the last layer use a ReLU activation function.

\subsection{Data}
Because the network is not linear, a much more robust data set is necessary in order to have the network discover a transformation that is sufficiently general. We trained neural networks on three different data sets. In all cases, the training data consists of 120,000 trajectories from Burgers' equation, each with 51 equally spaced time steps with $\Delta t = 0.002$. The validation data has the same structure as the training data but with 30,000 trajectories.

The difference in the three training data sets is the diversity of initial conditions. In data set 1, all of the initial conditions are randomly generated white noise. Recall that this was sufficient to find a global transformation for the heat equation. In data set 2, half of the initial conditions are white noise and the other half are sine waves. The sine waves are of the form $u_0(x) = A \sin(\omega x + \phi)$ where $A \in (0,1)$ and $\phi \in (0,2\pi)$ are chosen from a latin hypercube sample, and the frequency $\omega$ is an integer from 1 to 10 sampled from a truncated geometric distribution. In data set 3, one-third of the trajectories have a white noise initial condition, one-third have a sine wave initial condition, and one-third have a square wave initial condition where the height, width, and center of the square wave all come from a latin hypercube sample.

\begin{figure}[t]
\centering
\includegraphics[width=0.8\textwidth]{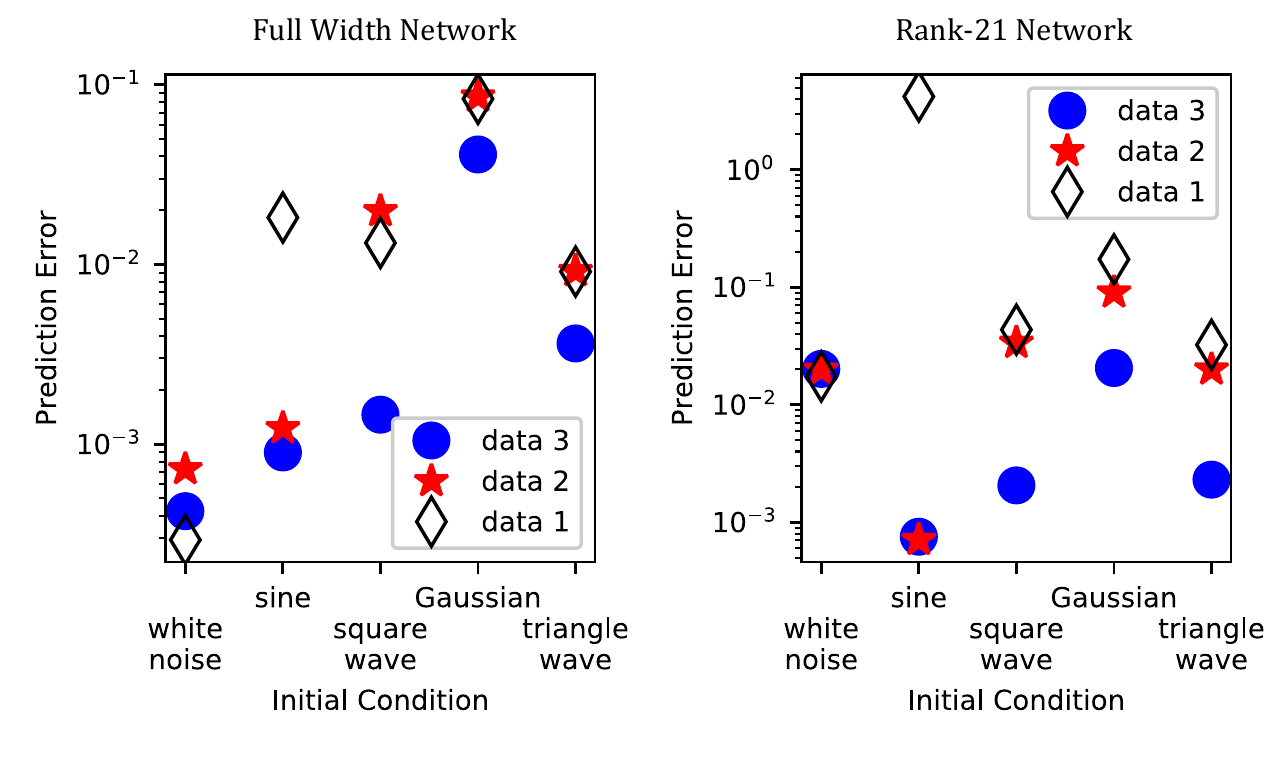}
 \caption{The network represented by blue circles was trained on white noise, sine, and square wave initial conditions. The red star network was trained on white noise and sines, and the black diamond network was just trained on white noise.}
 \label{fig:data_generalize}
\end{figure}

We trained networks on each of the three data sets. For each, we trained a full-width network with no dimensionality reduction in the middle layers, and we also trained a neural network with rank $r = 21$ in the middle layers. We then tested each network on test data that contains trajectories with five different types of initial conditions. The first three are the same as training data set 3 - white noise, sines, and square waves. The last two types of initial conditions are functions that were not seen in the training data for any of the data sets - Gaussians and triangle waves. The root-mean-square relative prediction errors for each network and each type of initial condition are shown in Figure \ref{fig:data_generalize}. The plot on the left is for full-width networks and the plot on the right is for the reduced rank $r = 21$. In both cases, the network trained on white noise has very low prediction error for trajectories that have white noise initial conditions, but does poorly on other types of initial conditions. Networks trained on white noise and sine waves have low prediction error for those two types of initial conditions, but still perform poorly on square wave initial conditions. The networks trained on white noise, sine waves, and square waves not only perform better than the other types of networks on predicting square wave initial conditions, they also generalize better to both Gaussians and triangle waves. This demonstrates that a variety of different initial conditions are needed in the training data in order to find a coordinate transformation that can represent any type of function.

\begin{figure}[t]
\centering
\includegraphics[width=0.7\textwidth]{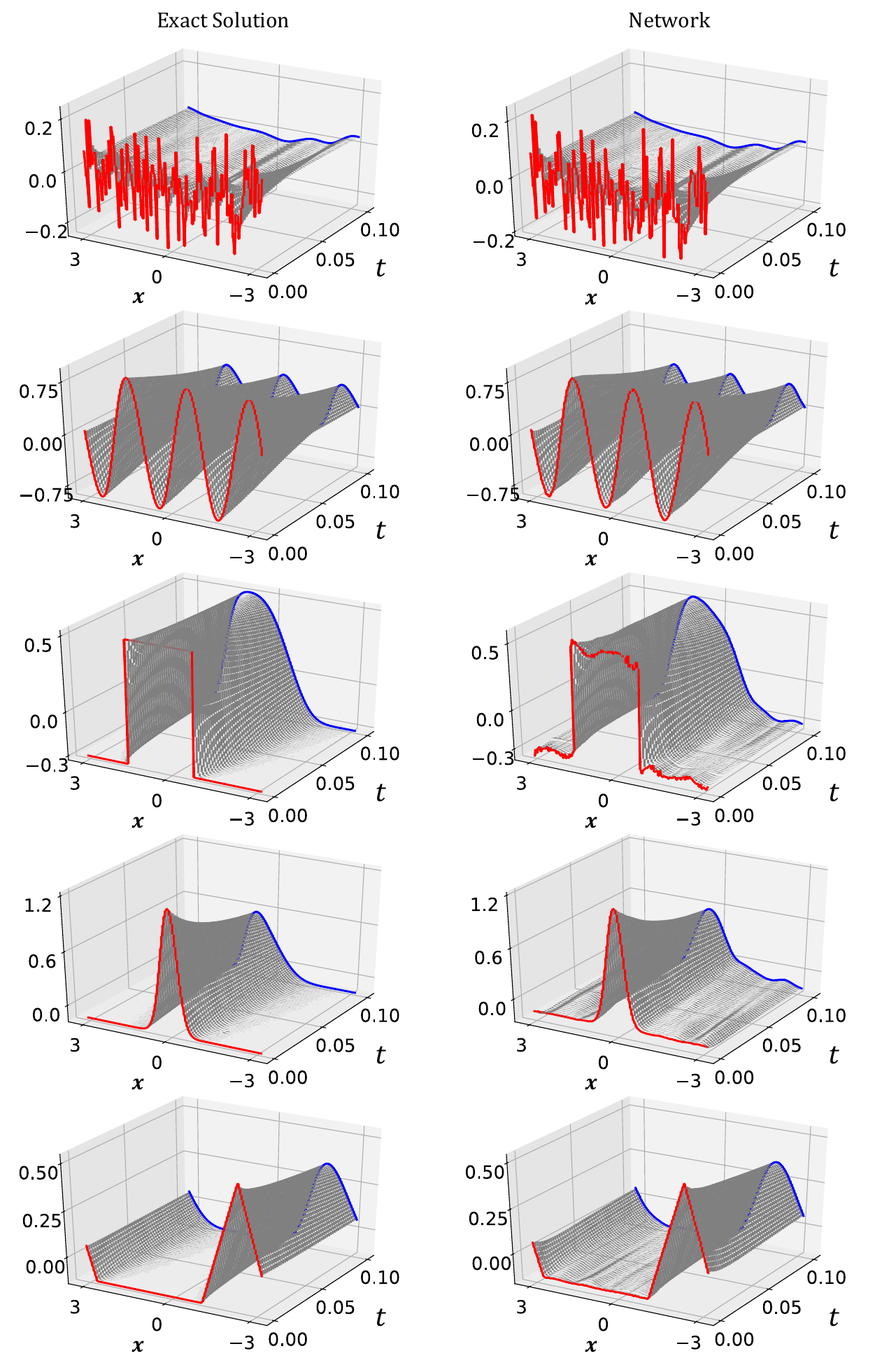}
 \caption{A comparison of the exact solution to Burgers' equation and the predictions given by the neural network. The top three use initial conditions of the same type as the training data, but the last two are types of initial conditions not found in the training data.}
 \label{fig:BurgersPred}
\end{figure}

Because the networks trained on data set 3 give better global transformations, all of the following results use that training and validation data. Figure \ref{fig:BurgersPred} shows the prediction given by the full-width network trained on data set 3 for a test trajectory with each of the five initial conditions in the test data. There is very good agreement between the exact solution and the network prediction in each case.

\subsection{Comparison with Cole-Hopf}
\begin{figure}[t]
\centering
\includegraphics[width=0.8\textwidth]{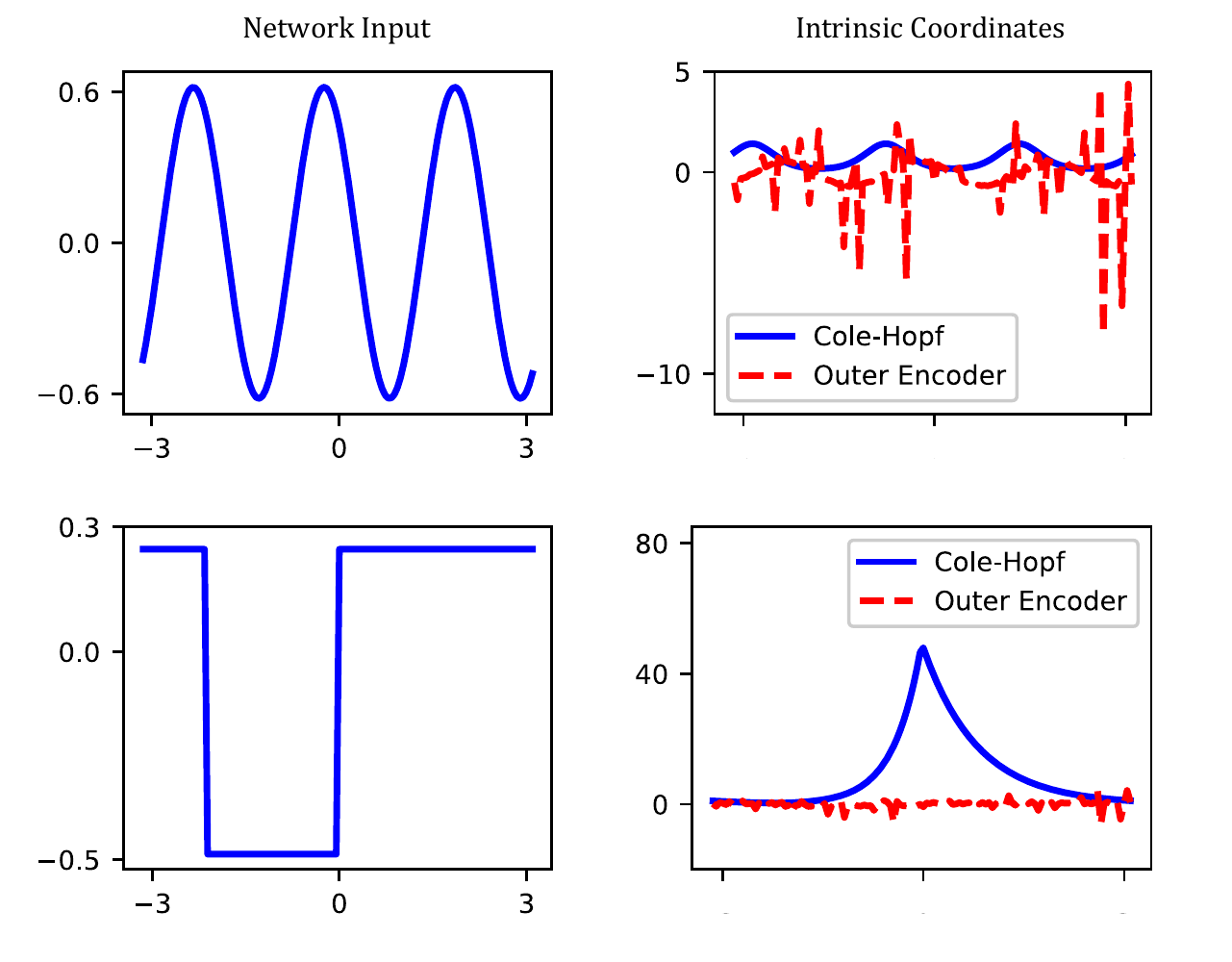}
 \caption{A comparison of the Cole-Hopf transformation given by equation \eqref{CHencode} and the output of the outer encoder of the neural network.}
 \label{fig:CHvsEncoded}
\end{figure}

Note that with $\epsilon = 10$ and $\mu = 1$, equation \eqref{CHencode} is $v(x,t) = e^{- 5 \int_0^x u(s,t) ds}$. Therefore, even if the Burgers' solution has an $\ell_1$ norm that is $\mathcal{O}(1)$, the function $v$ can have values on the order of $e^{5}$. This range of scales can be very difficult for a neural network to represent, especially with the use of regularization. However, similar to the heat equation, the Cole-Hopf transformation is not unique and can include scaling as long as that scaling is accounted for in the inverse transformation. Therefore, we do not expect the network to reproduce the formulas given by equations \eqref{CHencode} and \eqref{CHdecode}. Indeed, this is true. Figure \ref{fig:CHvsEncoded} shows a comparison between the Cole-Hopf transformation and the ``partially encoded'' function obtained by feeding the input through the network's outer encoder $\chi + \mathbf{I}$. The plots on the left show two functions $u(x,t)$. The plots on the right show the function $v(x,t)$ given by equation \eqref{CHencode} as a solid (blue) line and the output of the outer encoder as a dashed red line. The differences are stark. Notice for the square wave on the bottom that the Cole-Hopf transformation gives a function that is large in magnitude while the network outputs something much smaller in magnitude. However, the linear dynamics in this coordinate system still give good prediction for the Burgers' equation, as evidenced by Figure \eqref{fig:BurgersPred}.

\subsection{Reduced Order Model}
\begin{figure}[t]
\centering
\includegraphics[width=0.8\textwidth]{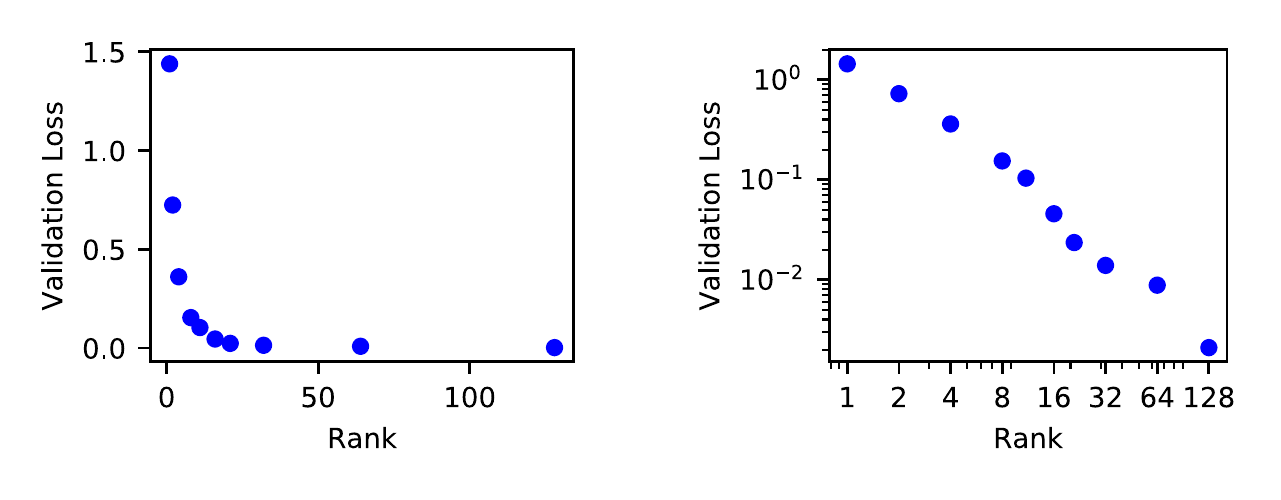}
 \caption{A plot of the validation losses for networks of different ranks.}
 \label{fig:ROM}
\end{figure}

By adjusting the widths of the middle layers of the network (labeled $\mathbf{v_k}$ and $\mathbf{v_{k+1}}$ in Figure \ref{fig:NetworkArch}), we can control the rank of the reduced order model of Burgers' equation. Recall that $r = 128$ is the full order model. We trained several networks with rank ranging from $r = 1$ to $r = 128$. Figure \ref{fig:ROM} shows the total validation loss for each network using a linear scale on the left and a logarithmic scale on the right. As expected, the loss decreases as the rank is increased.

\subsection{Poorly Selected Data or Architecture}

\begin{figure}[t]
\centering
\includegraphics{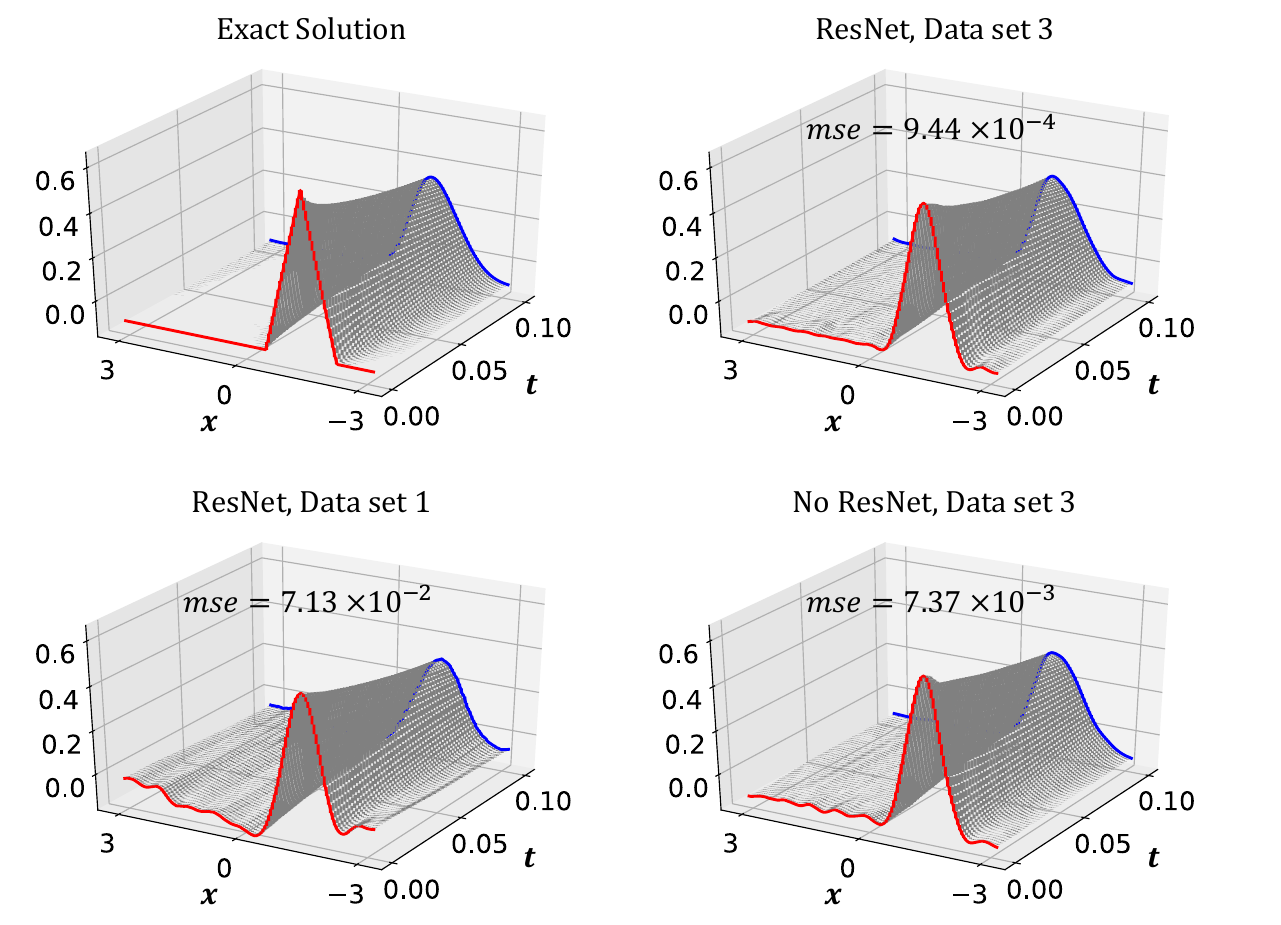}
 \caption{A comparison of the network predictions when the data is not sufficiently diverse, the architecture is not chosen properly, and when both the data and architecture are good.}
 \label{fig:BurgersCompare}
\end{figure}

When training a neural network for dynamics, it is important to use sufficiently diverse data as well as the proper architecture. Figure \ref{fig:BurgersCompare} shows the network predictions in the following cases. The top right plot shows the prediction from a neural network that is trained with the architecture given by Figure \ref{fig:NetworkArch} and data set 3. This gives the best result because the data has a diverse set of initial conditions and the architecture properly handles the identity. The bottom left plot shows the same architecture but trained on data set 1, which is less diverse than data set 3. The bottom right plot shows the prediction from a neural network with
the same architecture as Figure \ref{fig:NetworkArch} but without the skip connections that add the identity. All cases use a reduced order model with $r = 21$.

\section{Kuramoto-Sivashinsky Equation}

Our final example is the Kuramoto-Sivashinsky (KS) equation:
\begin{equation}
 u_t = -u u_x - u_{xx} - u_{xxxx}, \qquad x \in (-4 \pi,4\pi).
\end{equation}
As in the previous examples, we use periodic boundary conditions and a spatial discretization of 128 equally spaced points. Unlike Burgers' equation, there is no known coordinate transformation to linearize the KS equation. Therefore, the networks presented below provide the first known invertible transformation to linearize the KS equation.

\begin{figure}[t]
\centering
\includegraphics{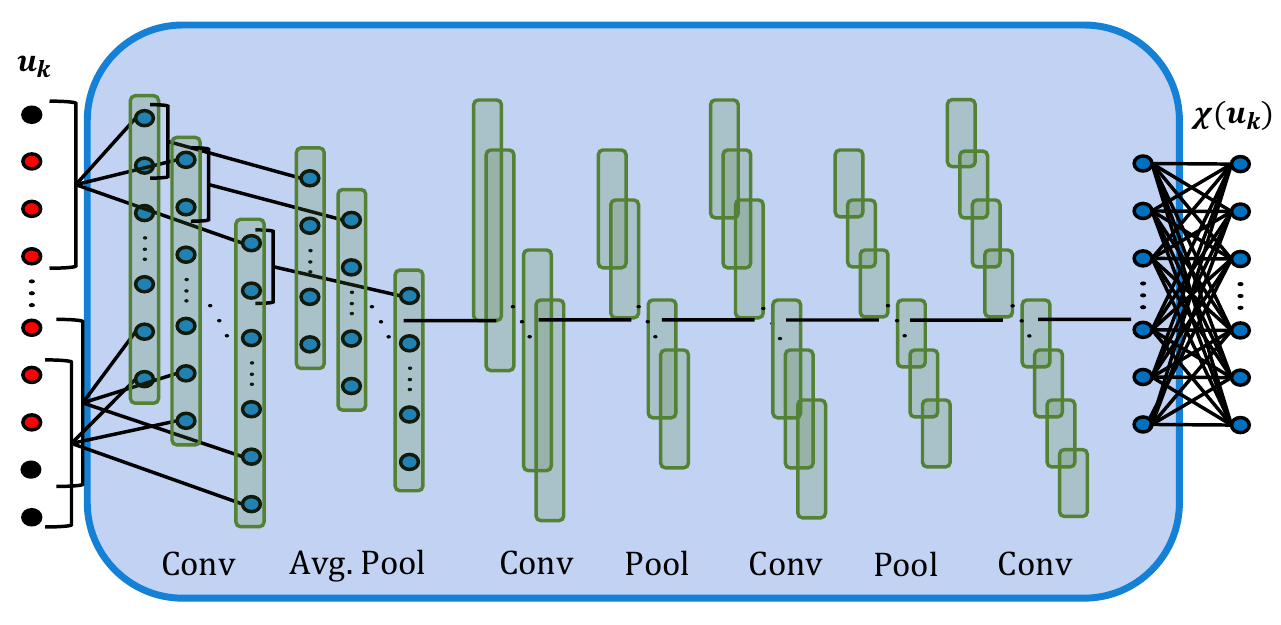}
 \caption{The network architecture for the outer encoder and decoder for the KS equation.}
\label{fig:KSArch}
\end{figure}

The neural network architecture is the same as that for Burgers' equation except for the architecture of the outer encoder and decoder. The architecture used for the KS equation is shown in Figure \ref{fig:KSArch}.  The first hidden layer is a convolutional layer containing 8 filters followed by an average pooling layer. This is then followed by a convolutional layer with 16 filters and an average pooling layer, a convolutional layer with 32 filters followed by an average pooling layer, and finally a convolutional layer with 64 filters. In all cases, the convolutional layers have kernel size 4, stride length 1, zero-padding, and ReLU activation functions while the average pooling layers have pool size 2, stride length 2, and no padding. The last convolutional layer is followed by a fully-connected layer with 128 neurons and ReLU activation and then a final fully-connected linear layer. 

The data set mirrors the data used for Burgers' equation. The training data consists of 120,000 trajectories, each with 51 equally spaced times steps. The initial conditions are evenly split between white noise, sine waves, and square waves. The validation data has the same structure as the training data but with 30,000 trajectories. The test data includes white noise, sine wave, and square wave initial conditions as well as Gaussian and triangle wave initial conditions.

We trained networks for several different data sets with the only difference being the time step. The time steps considered are $\Delta t = 0.125$, $\Delta t = 0.25$, $\Delta t = 0.5$, and $\Delta t = 1$. For each data set, we trained a full-width network (no dimensionality reduction) and a reduced-width network with rank $r = 21$. 

\begin{figure}[t]
\centering
\includegraphics[width=0.7\textwidth]{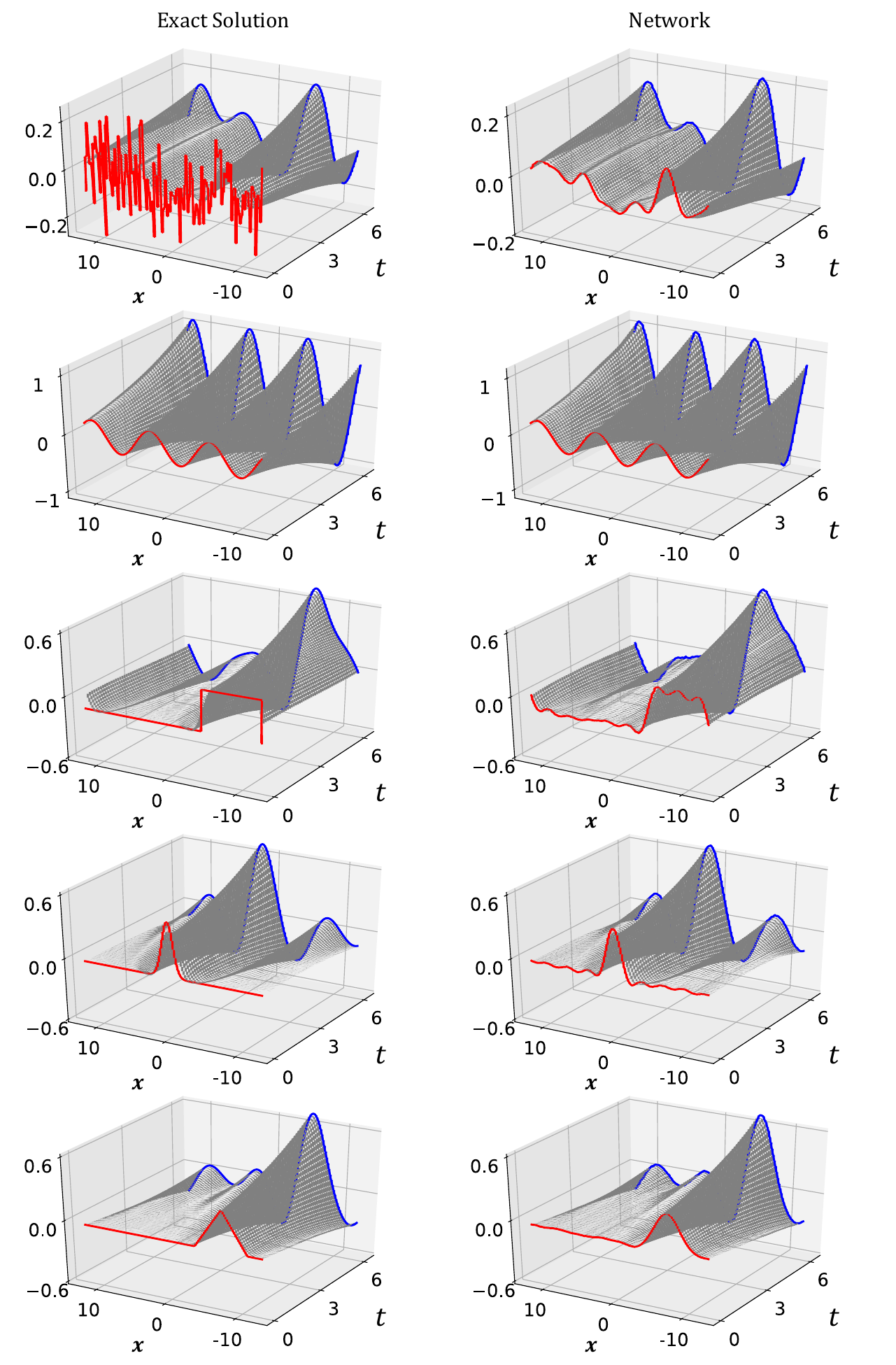}
 \caption{A comparison of the exact solution to the KS equation and the predictions given by the neural network. The top three use initial conditions of the same type as the training data, but the last two are types of initial conditions not found in the training data.}
 \label{fig:KSPred_smalldt}
\end{figure}

Figure \ref{fig:KSPred_smalldt} shows the predictions from the network trained with $\Delta t = 0.125$ and $r = 21$. The exact solution is shown on the left for five different initial conditions and the network output is on the right. The top three plots use initial conditions of the type used in the training data so good predictive power is expected. However, the bottom two plots use initial conditions that are not represented in the training data. The good agreement between the exact solution and the network predictions in these two cases show that the network gives a good global transformation to linearize the dynamics.

\begin{figure}[t]
\centering
\includegraphics[width=0.7\textwidth]{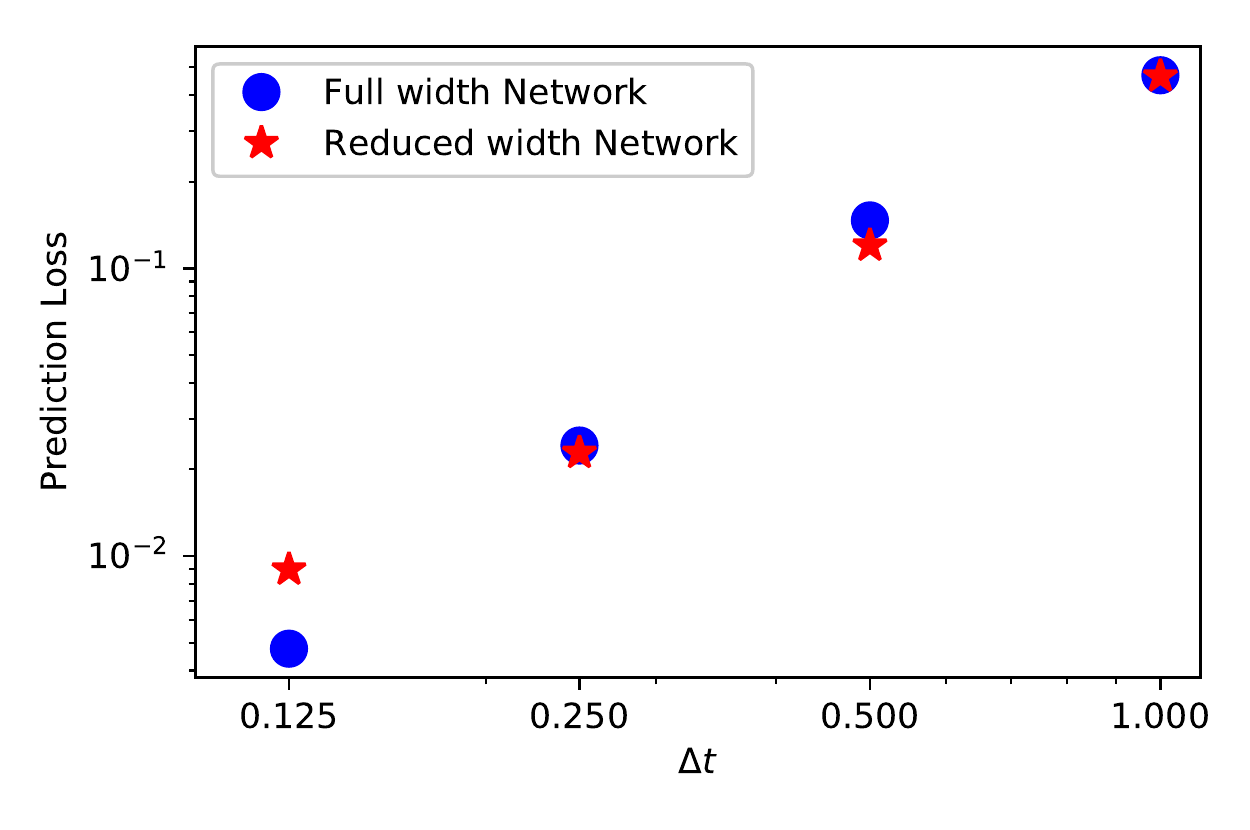}
 \caption{Prediction loss for the KS Equation for data with different time steps.}
 \label{fig:KS_timestep}
\end{figure}

Figure \ref{fig:KS_timestep} shows the prediction error on the test data for networks trained on data with different time steps. The results are given for both full width as well as reduced width networks. The prediction error is much lower for smaller time steps. This behavior is expected because a smaller time step leads to the solutions at consecutive times being closer together and therefore the network needs to approximate something nearer to the identity (and hence nearer to linear).

Note that the Cole-Hopf transformation does not depend on a time step. Similarly, the coordinate transformation for the KS equation should be independent of the time step, and therefore the encoder and decoder of the neural networks should be the same for data with any time step. The fact that smaller time steps are easier to represent with a neural network leads to a strategy for getting accurate predictions for larger time steps. In order to train the network for the data with $\Delta t = 0.25$, we initialized the weights with the weights from the network that had been trained with the shorter time step $\Delta t = 0.125$. This homotopy from the shorter time step to the longer time step led to both faster training and more accurate results. The total validation loss for $\Delta t = 0.25$ was $0.0192$ for the full width network and $0.0257$ for the reduced width network when using the homotopy strategy. When training the networks from scratch, the total validation losses were $0.0254$ for the full width network and $0.1753$ for the reduced width network. Therefore, failure to use a homotopy from the shorter time step led to a 32\% increase in the total validation loss for the full width network and a 582\% increase for the reduced width network. For $\Delta t = 0.5$, we used a homotopy from the $\Delta t = 0.25$ time step. Failure to homotopy led to a 45\% increase in the validation loss for the full width network and a 25\% increase for the reduced width network.




\section{Conclusion}
In this work, we have developed a scalable deep learning architecture specifically designed to learn a change of coordinates in which a given partial differential equation becomes linear.  
In particular, a custom deep autoencoder network is designed with additional constraints that in the low-dimensional latent space the system must evolve linearly in time.  
The resulting coordinate transformation is able to identify linearizing transformations that are analogous to known transformations, such as the Cole-Hopf transformation that maps the nonlinear Burgers equation into the linear heat equation. 
However the proposed architecture is designed to generalize to more complex PDE systems, such as the demonstrated Kuramoto-Sivashinsky.
Because our network is constrained to learn a linear dynamical system, the resulting network is highly interpretable, with the encoder network identifying eigenfunctions of the Koopman operator.
 
There are a number of future directions that arise from this work.  
The analysis performed here is promising and motivates the application of this approach on new PDE models for which known linearizing transformations do not exist.  
Extending these methods to higher-dimensional problems in two and three dimensions would also be interesting. 
Complex systems in fluid mechanics and turbulence, that exhibit multiscale phenomena in space and time, are especially interesting.  
These systems are generally characterized by a continuous frequency spectrum, motivating parameterized latent-space dynamics or time delays, which have both been successful for modeling ODEs with continuous spectra. 
Incorporating additional invariants and symmetries in the network architecture is another promising avenue of future work, as there is evidence that encoding partially known physics improves the learning rates and generalizability of the resulting models.  
The ability to embed nonlinear systems in a linear framework is particularly useful for estimation and control, where a wealth of techniques exist for linear systems.  
Therefore, it will likely be fruitful to extend these approaches to include inputs and control.

\begin{acknowledgments}
We would first like to thank Kathleen Champion, Jean-Christophe Loiseau, Karthik Duraisamy, Frank Noe, Yannis Kevrekidis and Igor Mezic for valuable discussions about sparse dynamical systems and autoencoders.
The authors acknowledge support from the Defense Advanced Research Projects Agency (DARPA PA-18-01-FP-125).  SLB further acknowledges support from the Army Research Office (ARO W911NF-17-1-0306), and JNK from the Air Force Office of Scientific Research (FA9550-17-1-0329). This work was facilitated through the use of advanced computational, storage, and networking infrastructure provided by AWS cloud computing credits funded by the STF at the University of Washington. This research was funded in part by the Argonne Leadership Computing Facility, which is a DOE Office of Science User Facility supported under Contract DE-AC02-06CH11357. 
\end{acknowledgments}


\bibliographystyle{siamplain}
\bibliography{references}

\end{document}